\documentclass[11pt]{amsart}
\usepackage{amsfonts,latexsym,amsthm,amssymb,graphicx}
\usepackage[all]{xy}

\DeclareFontFamily{OT1}{rsfs}{}
\DeclareFontShape{OT1}{rsfs}{n}{it}{<-> rsfs10}{}
\DeclareMathAlphabet{\mathscr}{OT1}{rsfs}{n}{it}

\def\demo{\smallskip\goodbreak{\it D\'emonstration.~~--~\kern.3em}
     \ignorespaces}%
\def\qedbox{$\square$}%
\def\qed{\ifmmode\qedbox\else\unskip\ \hglue0mm\hfill
     \qedbox\smallskip\goodbreak\fi}%
\def\proof{\smallskip\goodbreak{\it Proof\/.~~--~\kern.3em}
     \ignorespaces}%

\usepackage[english,francais]{babel}

\newtheorem{theoreme}{Th\'eor\`eme}[section]
\newtheorem{lemme}[theoreme]{Lemme}

\newtheorem{definition}[theoreme]{D\'efinition\rm}

\def\og{\leavevmode\raise.3ex\hbox{$\scriptscriptstyle\langle\!\langle$~}}
\def\fg{\leavevmode\raise.3ex\hbox{~$\!\scriptscriptstyle\,\rangle\!\rangle$}}

\title{
Classes de Chern des vari\'et\'es singuli\`eres, revisit\'ees
}
\author{Paolo Aluffi}
\address{Dept.~of Mathematics, Florida State University, Tallahassee
FL 32306, U.S.A.}
\address{Max-Planck-Institut f\"ur Mathematik, Postfach 7280,
D-53072 Bonn, Deutschland}
\email{aluffi@math.fsu.edu}

\begin{document}

\begin{abstract}
\selectlanguage{francais}
Nous introduisons une notion de groupe de proChow des vari\'et\'es,
qui reproduit la notion de groupe de Chow pour vari\'et\'es 
compl\`etes, et est fonctorielle par rapport aux morphismes arbitraires.
Nous construisons une transformation naturelle du foncteur des
fonctions constructibles vers le foncteur proChow, qui \'etend
la transformation naturelle de MacPherson. Nous illustrons
le r\'esultat en donnant des d\'emonstrations tr\`es courtes pour
(des g\'en\'eralisations de) deux faits bien connus sur le classes
de Chern-Schwartz-MacPherson.
\end{abstract}

\maketitle

\newcommand{\one}{1\hskip-3.5pt1}
\newcommand{\csm}{c_{\text{\rm SM}}}
\newcommand{\id}{{\text{id}}}
\newcommand{\A}{{\sf A}}
\newcommand{\F}{{\sf F}}
\newcommand{\T}{{\sf T}}
\newcommand{\cS}{{\mathscr S}}
\newcommand{\hA}{{\widehat {\sf A}}}
\newcommand{\oB}{{\overline B}}
\newcommand{\oU}{{\overline U}}
\newcommand{\oV}{{\overline V}}
\newcommand{\Zbb}{{\mathbb{Z}}}
\newcommand{\bxd}[1]{\{{#1}\}}

\section{Introduction}
Soit $X$ une vari\'et\'e sur un corps alg\'ebriquement clos de 
caract\'eristique nulle.  Des notions \'equivalentes de 
{\em classe de Chern totale\/} de $X$ ont \'et\'e donn\'ees 
ind\'ependamment par Marie-H\'el\`ene Schwartz (\cite{MR32:1727})
et par Robert MacPherson (\cite{MR50:13587}) 
pour les vari\'et\'es alg\'ebriques complexes compactes, en homologie;  
la d\'efinition a \'et\'e depuis \'etendue aux vari\'et\'es alg\'ebriques 
compl\`etes sur tous corps alg\'ebriquement clos de caract\'eristique 
nulle, dans le groupe de Chow $\A_*X$ de $X$.  
Nous appelons cette classe la {\em classe de 
Chern-Schwartz-MacPherson\/} (CSM) de $X$, d\'enot\'ee $\csm(X)$.

La classe de CSM co\"incide avec la classe de Chern ordinaire
du fibr\'e tangent quand $X$ est lisse: $\csm(X)=c(TX)\cap [X]$ 
dans ce cas. Elle satisfait \'egalement une propri\'et\'e remarquable 
de fonctorialit\'e: elle est d\'efinie comme la valeur $c_*(\one_X)$ 
prise par une transformation {\em naturelle\/} $c_*$ du foncteur 
des fonctions constructibles $\F$ vers le foncteur de Chow $\A_*$,
sur la fonction constante $\one_X$.  Alexandre Grothendieck et 
Pierre Deligne avaient conjectur\'e l'existence de cette 
transformation naturelle;  MacPherson la construit explicitement 
dans \cite{MR50:13587}, en utilisant d'autres invariants importants
qu'il introduit dans ce travail. 

Dans cette note nous proposons une construction alternative
des classes de Chern-Schwartz-MacPherson, dans un groupe `enrichi'
de Chow $\hA_*X$ (le groupe de {\em proChow\/}) obtenu en prenant 
des limites appropri\'ees des groupes de Chow ordinaires.  Le groupe 
de proChow est un foncteur covariant par rapport \`a tous les 
morphismes (tandis que le groupe de Chow ordinaire est fonctoriel
seulement par rapport aux morphismes {\em propres\/}); nous 
d\'emontrons que la {\em transformation de proCSM\/} correspondante 
$\F \leadsto \hA_*$ est naturelle par rapport aux morphismes
quelconques. Si $X$ est compl\`ete, le groupe de proChow de $X$ est 
canoniquement isomorphe au groupe de Chow ordinaire, et sa {\em
classe proCSM\/} est \'egale \`a la classe de
Chern-Schwartz-MacPherson.

Notre d\'efinition est directe, sans r\'ef\'erence aux invariants 
auxiliaires comme les classes de Chern-Mather ou l'obstruction 
locale d'Euler. Pour illustrer son utilisation nous donnons 
des d\'emontrations tr\`es condens\'ees de deux r\'esultats connus sur 
les classes de CSM, parus dans les {\em C.R.A.S.\/}:  la {\em formule 
du produit\/} de Kwieci\'nski (\cite{MR1158750}), et la formule
de Ehlers-Barthel-Brasselet-Fieseler pour les classes de CSM des 
vari\'et\'es toriques (\cite{MR1197235}).

\newpage
\section{Foncteur proChow}\label{proChow}

Nous travaillons sur un corps alg\'ebriquement clos $k$ de 
caract\'eristique nulle.  

Soit $\cS$ une cat\'egorie de $k$-vari\'et\'es.  Pour $U$ dans $\cS$, 
soit $\cS_U$ la cat\'egorie dont les objets sont les $\cS$-morphismes 
$i: U \to Z^i$ de $U$ vers les vari\'et\'es {\em compl\`etes\/} dans $\cS$, 
et les morphismes $j\mapsto i$ sont les diagrammes commutatifs 
de $k$-vari\'et\'es
$$\xymatrix@R=0pt{
& {Z^j} \ar[dd]^\pi \\
U \ar[ur]^j \ar[dr]_i \\
& {Z^i}
}$$
o\`u $\pi$ est un morphisme {\em propre.\/} Nous supposons que les 
conditions suivantes sur $\cS$ sont v\'erifi\'ees:  
\begin{itemize}
\item pour tout $U$ dans $\cS$ et toute paire d'objets $i,j$ de $\cS_U$, 
il y a un objet $k$ dans $\cS_U$ tel que $k \to i$ et $k \to j$, et $k: U
\hookrightarrow Z^k$ est une {\em adh\'erence\/} (c'est-\`a-dire, 
$k$ est un plongement ouvert, et $\oU=Z^k$);  
\item si $U$ est lisse, on peut choisir l'adh\'erence $Z^k$ 
comme ci-dessus et {\em bonne:\/} c'est-\`a-dire, $Z^k$ est lisse,
et le compl\'ement $Z^k\smallsetminus U$ est un diviseur \`a 
croisements normaux et composantes lisses.
\end{itemize}

Par exemple, ces conditions sont satisfaites pour la cat\'egorie 
de {\em toutes\/} les $k$-vari\'et\'es (en caract\'eristique nulle, par 
r\'esolution des singularit\'es).  

\begin{definition} 
Le {\em groupe de proChow\/} de $U$ (par rapport a $\cS$) est la
limite $\hA^\cS_*U:=\varprojlim_i \A_* Z^i$.
\end{definition}

Concr\`etement, un \'el\'ement $\rho\in \hA^\cS_* U$ consiste en le 
choix d'un \'el\'ement $\rho^i$ dans le groupe de Chow (conventionnel) 
$\A_*Z^i$ pour tout $i$ dans $\cS_U$, sujet \`a la condition de 
compatibilit\'e $\pi_*\rho^j=\rho^i$ pour tout $\pi: j \to i$.  
Nous disons que $\rho^i$ est la {\em composante de $\rho$ dans $\A_*Z^i$.\/}

Nous omettrons l'indice sup\'erieur $\cS$ quand aucune ambigu\"it\'e 
n'est probable; le lecteur notera que le groupe de proChow d\'epend 
de la cat\'egorie choisie $\cS$.  Les faits suivants sont cependant
ind\'ependants de $\cS$ (si $\cS$ satisfait aux conditions de 
`cofinalit\'e' sp\'ecifi\'ees ci-dessus) et imm\'ediatement v\'erifi\'es:

\begin{lemme}\label{simpl}
Avec les notations ci-dessus:
\begin{itemize}
\item Si $U$ est compl\`ete, il y a un isomorphisme canonique 
$\hA_* U \cong \A_*U$.
\item Pour sp\'ecifier un \'el\'ement de $\hA_*U$, il suffit de
choisir un ensemble compatible de $\rho^i\in \A_*Z^i$ pour 
toute adh\'erence $i:U \to Z^i$ dans $\cS$.
\item Si, en outre, $U$ est lisse, il suffit de choisir 
un ensemble compatible de $\rho^i\in \A_*Z^i$ pour toute
{\em bonne\/} adh\'erence $i:U \to Z^i$ dans $\cS$.
\end{itemize}
\end{lemme}

Chaque sous-sch\'ema $B$ de $U$ d\'etermine un \'el\'ement distingu\'e 
$[\overline B]$ de $\hA_*U$: pour chaque adh\'erence $j: U \to Z^j$, 
on choisit la classe $[\overline B] \in \A_*Z^j$ de l'adh\'erence de
$B$ dans $Z^j$; ce choix est clairement compatible. Si $U$ est 
compl\`ete, $[\oB]\in\hA_*U\cong \A_*U$ est la `classe fondamentale'
ordinaire de $\oB$.

Le groupe de proChow $\hA_*=\hA_*^\cS$ est un foncteur $\cS\leadsto $
Groupes Ab\'eliens: si $f: X \to Y$ est un morphisme dans $\cS$, 
alors $j \to j\circ f$ induit un foncteur $\cS_Y \to \cS_X$, et 
donc un homomorphisme $f_*:\hA_*X \to \hA_*Y$. Concr\`etement, pour 
$\rho\in \hA_*X$ et $j:Y \to Z^j$ dans $\cS_Y$, la composante de
$f_* \rho$ dans $\A_* Z^j$ est simplement \'egale \`a la composante de 
$\rho$.  Si $f$ est propre et $X$ et $Y$ sont compl\`etes, 
alors $f_*: \hA_* X\cong \A_*X \to \A_*Y\cong \hA_*Y$ est le
push-forward propre ordinaire des groupes de Chow. On note
cependant que tandis que $\A_*$ est fonctoriel seulement par
rapport aux morphismes propres, le proChow $\hA_*$ est fonctoriel
par rapport \`a tous les morphismes dans $\cS$.

\section{Classes proCSM}\label{proCSM}

Avec $\cS$ comme dans la Section~\ref{proChow}, et $X$ dans $\cS$, nous 
d\'efinissons le groupe des {\em fonctions $\cS$-constructibles\/} 
$\F^\cS(X)$ comme le groupe des combinaisons finies $\Zbb$-lin\'eaires 
des fonctions caract\'eristiques $\one_U$ (o\`u $\one_U(p)=1$ si $p\in U$, 
et $0$ si $p\in X\smallsetminus U$) o\`u $U$ sont des sous-vari\'et\'es 
localement ferm\'ees {\em lisses\/} de $X$, telles que les
inclusions $U\subset X$ sont des morphismes de $\cS$.

Nous posons maintenant une condition suppl\'ementaire sur $\cS$.  
Nous demandons que le push-forward conventionnel des fonctions 
constructibles (d\'efinies en prenant la caract\'eristique d'Euler 
fibre \`a fibre, voir \cite{MR50:13587} pour le cas complexe) 
pr\'eserve la $\cS$-constructibilit\'e: c'est-\`a-dire, qu'il d\'efinisse
un push-forward $f_*: \F^\cS(X) \to \F^\cS(Y)$ pour chaque morphisme 
$f: X\to Y$ dans $\cS$. Nous exigeons \'egalement que $\one_X$ soit 
$\cS$-constructible pour chaque $X$ dans $\cS$.

Sous ces conditions, $\F^\cS$ d\'efinit (en caract\'eristique nulle!)  
un foncteur covariant $\cS \leadsto$ Groupes Ab\'eliens, et chaque 
$X$ dans $\cS$ d\'etermine un \'el\'ement distingu\'e $\one_X\in \F^\cS(X)$.  
Nous omettrons habituellement l'indice sup\'erieur $\cS$.

Nous d\'efinissons maintenant un homomorphisme $\F(X) \to 
\hA_*X$, $\alpha \mapsto \bxd{\alpha}$, et un \'el\'ement distingu\'e 
$\bxd X:=\bxd {\one_X}\in \hA_*X$. Nous commen\c cons par le cas 
lisse:

\begin{definition}\label{defnonsing}
Soit $U$ lisse, dans $\cS$. La {\em classe 
proCSM\/} de $U$ dans $\hA_*U$, not\'ee $\bxd U$, est l'\'el\'ement 
du groupe de proChow d\'etermin\'e par 
$c(\Omega^1_{\oU}(\log D)^\vee)\cap [\oU]\in \A_*\oU$
pour toute bonne adh\'erence $\overline U$ de $U$ dans~$\cS_U$, 
o\`u $D=\oU\smallsetminus U$ est le diviseur \`a croisements normaux 
correspondant, et $\Omega^1_\oU(\log D)^\vee$ d\'esigne le 
dual du fibr\'e de formes diff\'erentielles avec p\^oles logarithmiques 
le long de~$D$.  
\end{definition}

Ce choix est compatible selon la section~\ref{proChow}, comme on le
d\'emontrera dans le Th\'eor\`eme~\ref{main}; donc il d\'efinit
un \'el\'ement de $\hA_*U$, d'apr\`es le Lemme~\ref{simpl}.

Soit maintenant $X$ une vari\'et\'e arbitraire (c'est-\'a-dire,
pas n\'ecessairement lisse) dans $\cS$, et soit $\alpha\in F(X)$
une fonction constructible sur $X$. Soit $\alpha=\sum_U m_U
\one_U$, avec $U$ lisse, localement ferm\'ee, $i_U: U\subset X$ 
dans $\cS$, et $m_U\in \Zbb$.

\begin{definition}\label{defpossing}
La {\em classe proCSM de $\alpha$\/} est la somme
$\bxd{\alpha}=\sum_U m_U {i_U}_* \bxd{U}\in \hA_*X$.
La {\em classe proCSM de $X$\/} est la classe $\bxd{X}:=\bxd{\one_X}$.
\end{definition}

Nous pouvons maintenant \'enoncer et d\'emontrer le r\'esultat
principal de cette note.

\begin{theoreme}\label{main}
Avec ces notations:
\begin{enumerate}
\item\label{1} Les choix donn\'es dans la D\'efinition~\ref{defnonsing} 
sont compatibles: c'est-\`a-dire, si $i: U \to \oU^i$ et $j: U \to 
\oU^j$ sont des bonnes adh\'erences de $U$ dans $\cS_U$, de 
compl\'ements $D^i$, $D^j$, et $\pi: \oU^j \to \oU^i$ est un morphisme 
propre tels que $i=\pi\circ j$, alors
$\pi_* \left(c(\Omega^1_{\oU^j}(\log D^j)^\vee)\cap [\oU^j]\right)
=c(\Omega^1_{\oU^i}(\log D^i)^\vee)\cap [\oU^i]$.
\item\label{2} La D\'efinition~\ref{defpossing} est ind\'ependante
des choix: c'est-\`a-dire, si $\alpha=\sum_U m_U \one_U = \sum_V n_V 
\one_V$ sont deux mani\`eres d'exprimer $\alpha$ comme combinaison
lin\'eaire finie des fonctions caract\'eristiques des sous-vari\'et\'es
localement ferm\'ees lisses de  $X$, alors
$\sum_U m_U {i_U}_* \bxd{U} = \sum_V n_V {i_V}_* \bxd{V}$.
\item\label{3} L'homomorphisme $\F(X) \to \hA_*X$, $\alpha\mapsto 
\bxd{\alpha}$ donn\'e dans la D\'efinition~\ref{defpossing} d\'efinit
une transformation {\em naturelle\/} $\F \leadsto \hA_*$; c'est-\`a-dire,
$f_* \bxd{\alpha} = \bxd{ f_*(\alpha)}$
pour chaque morphisme $f:X \to Y$ in $\cS$.
\item\label{4} Si $X$ est {\em compl\`ete,\/} alors la classe proCSM
de $X$ est la classe de Chern-Schwartz-MacPherson conventionnelle:
$\bxd{X} = \csm(X) \in \A_*X \cong \hA_*X$.
\end{enumerate}
\end{theoreme}

On peut d\'emontrer le trois premiers points ind\'ependamment
du r\'esultat de MacPherson en \cite{MR50:13587}; cela est fait
dans \cite{math.AG/0507029}. Dans le cas particulier des 
vari\'et\'es compl\`etes et des fonctions propres, le troisi\`eme
point donne une transformation naturelle comme prescrite
par la (version Chow de la) conjecture de Grothendieck-Deligne;
l'\'egalit\'e des classes proCSM et des classes CSM pour les
vari\'et\'es compl\`etes d\'ecoule alors de l'unicit\'e de cette
transformation naturelle (qui est une cons\'equence imm\'ediate
de la r\'esolution des singularit\'es).

Nous pr\'esentons ici une preuve qui utilise le th\'eor\`eme de
MacPherson (c'est-\`a-dire, le fait que la transformation
de MacPherson $c_*$ est naturelle), ce qui permet un raisonnement
plus rapide.

\begin{demo}
Si $U$ est lisse, $\oU$ est une bonne adh\'erence de $U$, et
$D=\oU\smallsetminus U$, alors
\begin{equation*}
\tag{\dag}
c(\Omega^1_{\oU}(\log D)^\vee)\cap 
[\oU]=c_*(\one_U)\in A_*\oU
\quad.
\end{equation*}
Ceci d\'ecoule facilement du fait que $c_*$ est naturel et d'un
calcul explicite de classes de Chern; voir par exemple 
la Proposition~15.3 dans \cite{MR1893006}, ou le Th\'eor\`eme~1 
dans \cite{MR2001d:14008}.

({\it \ref{1}\/}) d\'ecoule de (\dag), du fait que $c_*$ est naturel,
et de la d\'efinition de push-forward des fonctions constructibles:
$$\pi_* \left(c(\Omega^1_{\oU^i}(\log D^i)^\vee)\cap [\oU^i]\right)
=\pi_* c_*(\one_U)=c_* \pi_*(\one_U) = c_*(\one_U)
=c(\Omega^1_{\oU^j}(\log D^j)^\vee)\cap [\oU^j]\quad.$$

La preuve des autres points est simplifi\'ee par la
version alternative suivante de la {D\'efini\-tion}~\ref{defpossing}:

\begin{lemme}\label{key}
Pour chaque $z:X \to Z$ dans $\cS_X$, et chaque $\alpha\in \F(X)$,  
la composante de $\bxd{\alpha}$ dans $\A_*Z$ est $c_* (z_*(\alpha))$.
\end{lemme}

Pour d\'emontrer le lemme, on utilise (\dag) pour \'ecrire la composante
de $\bxd U$ dans $\A_*Z$, pour $U$ lisse et $z_U:U \to Z$ dans $\cS$,
comme $c_*({z_U}_*(\one_U))$. Si $\alpha\in \F(X)$, alors 
$\alpha=\sum_U m_U \one_U$, avec $U$ lisse et l'inclusion
$U\subset X$ dans $\cS$; donc, pour chaque $z:X \to Z$, la composante
de $\sum_U m_U \bxd U$ dans $\A_*Z$ est $\sum_U m_U c_*(z_*\one_U)= 
c_*(z_*(\sum_U m_U \one_U))= c_*(z_*(\alpha))$, comme \'enonc\'e.

({\it \ref{2}\/}) s'ensuit imm\'ediatement, puisque $c_* (z_*(\alpha))$
d\'epend seulement de $\alpha$, et pas de la d\'ecomposition
$\alpha=\sum_* m_U \one_U$. D'ailleurs, si $X$ est elle-m\^eme
compl\`ete, alors en prenant $Z=X$, $z=\id_X$, et
$\alpha=\one_X$ dans le Lemme~\ref{key}, on trouve ({\it \ref{4}\/}).

Enfin, le Lemme~\ref{key} implique ({\it \ref{3}\/}).
En effet, soit $z:Y \to Z$ un objet quelconque de $\cS_Y$; alors 
$w\circ f$ est un objet de $\cS_X$ et, d'apr\`es la d\'efinition
de push-forward de groupes de proChow, la composante dans $\A_*Z$ 
du push-forward $f_* \bxd{\alpha}$ est simplement \'egale \`a la
composante de $\bxd{\alpha}$ dans $\A_*Z$. D'apr\`es le Lemme~\ref{key}, 
cette composante est $c_*((w\circ f)_* \alpha)=c_*(z_* (f_*(\alpha)))$,
et encore d'apr\`es le Lemme~\ref{key}  elle est \'egale \`a la composante
de $\bxd{f_*(\alpha)}$ dans $\A_*Z$, ce qui prouve ({\it \ref{3}\/}).
\end{demo}

\section{Exemples}\label{examples}

Nous illustrons le formalisme pr\'esent\'e dans \S\ref{proChow} et 
\S\ref{proCSM} en donnant des preuves tr\`es courtes (et valides dans le 
cadre proCSM) de deux r\'esultats connus sur les classes de 
Chern-Schwartz-MacPherson.

Nous utiliserons des cat\'egories diff\'erentes, comme permis par les
constructions donn\'ees dans les sections pr\'ec\'edentes. En notant par 
$\hA_*$ le foncteur proChow obtenu \`a partir de $\cS=$ la cat\'egorie de 
toutes $k$-vari\'et\'es, et par $\F$ le foncteur de fonctions constructibles
sur cette cat\'egorie, on note que pour toute autre cat\'egorie $\cS$ on 
a des homomorphismes canoniques $\F^\cS(X) \to \F(X)$, $\hA_*(X) \to
\hA^\cS_*(X)$, compatibles avec les transformations naturelles proCSM 
correspondantes.

Le premier r\'esultat est la {\em formule de produit\/} de Micha{\l } Kwieci{\'n}ski.
On prend pour $\cS$ la cat\'egorie de {\em produits\/} $X\times Y$ 
(techniquement, de paires $(X,Y)$), o\`u $X$ et $Y$ sont des $k$-vari\'et\'es, et 
dont les morphismes $X_1 \times Y_1 \rightarrow X_2 \times Y_2$ consistent 
en des paires $(f, g)$, o\`u $f:X_1\to X_2$ et $g:Y_1 \to Y_2$ sont morphismes.
Les conditions d\'etaill\'ees dans \S2 et \S3 sont clairement v\'erifi\'ees 
(en caract\'eristique nulle), et on a donc un foncteur $\hA^\times_*$ de proChow 
et un foncteur $\F^\times$ de fonctions $\cS$-constructibles.  Le groupe 
$\F^\times(X\times Y)$ consiste en des fonctions $\alpha\otimes \beta$ d\'efinies 
par $\alpha\otimes\beta(x,y)=\alpha(x)\beta(y)$, o\`u $\alpha \in \F(X)$, $\beta
\in \F(Y)$ sont des fonctions constructibles (ordinaires). On note la classe 
proCSM correspondante par $\bxd{\alpha \otimes\beta}^\times$.

On voit en outre qu'il y a un homomorphisme canonique \'evident
$$\xymatrix{
\hA_*(X) \otimes \hA_*(Y) \ar[r]^\otimes &  \hA^\times_*(X\times Y)
}\quad,$$
$(\alpha,\beta) \mapsto \alpha \otimes \beta$,
induit par les produits ext\'erieurs pour les groupes de Chow ordinaires 
(\cite{MR85k:14004}, \S1.10).

\begin{theoreme}
Soient $X$ et $Y$ deux vari\'et\'es, $\alpha\in \F(X)$, $\beta\in \F(Y)$
et $\alpha\otimes\beta \in \F^\times (X\times Y)$ comme ci-dessus. Alors
$\bxd{\alpha \otimes \beta}^\times = \bxd \alpha \otimes \bxd \beta$.
\end{theoreme}

\begin{demo}
Par bilin\'earit\'e, l'\'enonc\'e suit du cas $\alpha=\one_U$, $\beta=\one_V$ 
pour $U$,~$V$ sous-vari\'et\'es lisses de $X$, $Y$ resp.; c'est-\`a-dire qu'il suffit
de v\'erifier que pour $U$, $V$ lisses, et pour $\oU$, $\oV$ bonnes adh\'erences
of $U$, $V$, de compl\'ements $D=\oU \smallsetminus U$,
$E=\oV\smallsetminus V$,
$$c(\Omega^1_{\oU\times \oV}(\log(D+E))^\vee)\cap [\oU\times \oV]
=\left(c(\Omega^1_\oU(\log D)^\vee)\cap [\oU]\right) \otimes
\left(c(\Omega^1_\oV(\log E)^\vee)\cap [\oV]\right)
\quad,$$
et ceci suit imm\'ediatement du calcul standard des classes de Chern pour
le fibr\'e de formes diff\'erentielles avec p\^oles logarithmiques. 
\end{demo}

Le cas particulier dans lequel $X$ et $Y$ sont compl\`etes reproduit le 
th\'eor\`eme de Kwieci{\'n}ski (\cite{MR1158750}), puisque dans ce cas-l\`a 
tous les groupes de proChow dans l'\'enonc\'e sont isomorphes aux groupes 
conventionnels de Chow (d'apr\`es le Lemme~\ref{simpl}), et les classes 
proCSM sont \'egales aux classes de Chern-Schwartz-MacPherson 
(d'apr\`es le Th\'eor\`eme~\ref{main}).

Notre deuxi\`eme exemple est la formule de Fritz Ehlers pour la classe de 
Chern-Schwartz-MacPherson d'une vari\'et\'e torique;  voir \cite{MR1234037}, 
p.~113, et \cite{MR1197235} pour la preuve dans le cadre CSM conventionnel. 
Pour donner un \'enonc\'e et une preuve dans le cas plus g\'en\'eral de proChow, 
soit $\cS$ la cat\'egorie de $k$-vari\'et\'es toriques, avec les morphismes 
$T$-\'equivariants.
Le foncteur et les classes proCSM correspondants seront not\'es 
respectivement $\hA^\T$ et $\bxd{X}^\T$; la {\em classe fondamentale\/}
de $B\subset X$ dans $\hA^\T_*(X)$ sera not\'ee $[\oB]$.

\begin{theoreme}\label{ehlers}
Soit $X$ une vari\'et\'e torique. Alors 
$\bxd{X}^\T=\sum_{B\in X/T} [\oB]\in \hA^\T_*(X)$,
o\`u la somme porte sur l'ensemble (fini) des $T$-orbites.
\end{theoreme}

\begin{demo}
Puisque $X$ est l'union des $T$-orbites $B$ nous avons $\bxd X^\T=
\sum\bxd B^\T$, et par cons\'equent il suffit de montrer que si $B$ est l'orbite 
ouverte dans la sous-vari\'et\'e torique $\overline B\subset X$ alors
$\bxd B^\T=[\oB] \in \hA_*^\T(B)$;
ce qui est \'equivalent \`a montrer que si $\oB$ est une { \em bonne\/} 
adh\'erence (torique) de $B$, et $D=\oB\smallsetminus B$,
$c(\Omega^1_{\overline B}(\log D)^\vee)\cap [\overline B]=[\overline B]
\in \A_*(\overline B)$:
et ceci est vrai puisque $\Omega^1_{\overline B}(\log D)$ est trivial 
(\cite{MR1234037}, Proposition, p.~87).
\end{demo}

Dans le cas particulier o\`u $X$ est une vari\'et\'e torique {\em compl\`ete,\/} 
on retrouve la formule d'Ehlers.

Le th\'eor\`eme~\ref{ehlers} admet (avec la m\^eme preuve) une 
g\'en\'eralisation aux plongements toro\"idaux qui ne sont pas 
n\'ecessairement normaux.

\end{document}